\documentclass{article}

\usepackage[a4paper]{geometry}

\usepackage[francais,english]{babel}
\usepackage{amsfonts,latexsym,amstext}
\usepackage{amsmath,amscd,amssymb}

\usepackage{enumerate}

\usepackage[latin1]{inputenc}

\usepackage{aeguill}
\usepackage[T1]{fontenc}

\usepackage[usenames]{color}

\usepackage[a4paper]{geometry}

\usepackage{url}

\usepackage{psfrag}

\usepackage{placeins}

\usepackage[dvips]{graphicx}
\usepackage{latexsym,amssymb,amsmath,color,times,bbm}
\newtheorem{thm}{Theorem}[section]

\newtheorem{lem}[thm]{Lemma}
\newtheorem{prop}[thm]{Proposition}
\newtheorem{rem}[thm]{Remark}

\numberwithin{equation}{section}
\newcommand{\norm}[1]{\left\Vert#1\right\Vert}
\newcommand{\abs}[1]{\left\vert#1\right\vert}
\newcommand{\set}[1]{\left\{#1\right\}}

\def\R{\mathbb R}
\def\N{\mathbb N}

\def\E{\mathbb E}
\def\P{\mathbb P}
\def\Q{\mathbb Q}

\def\cprime{$'$}

\newcommand{\tr}[1]{{\vphantom{#1}}^{\mathit t}{#1}}

\def\sqw{\hbox{\rlap{\leavevmode\raise.3ex\hbox{$\sqcap$}}$%
\sqcup$}}
\def\sqb{\hbox{\hskip5pt\vrule width4pt height6pt depth1.5pt%
\hskip1pt}}

\def\qed{\ifmmode\hbox{\hfill\sqb}\else{\ifhmode\unskip\fi%
\nobreak\hfil
\penalty50\hskip1em\null\nobreak\hfil\sqb
\parfillskip=0pt\finalhyphendemerits=0\endgraf}\fi}
\def\cqfd{\ifmmode\sqw\else{\ifhmode\unskip\fi\nobreak\hfil
\penalty50\hskip1em\null\nobreak\hfil\sqw
\parfillskip=0pt\finalhyphendemerits=0\endgraf}\fi}
\begin{document}

\renewcommand{\labelitemi}{$\bullet$}
\bibliographystyle{plain}
\pagestyle{headings}
\title{A note on the existence of solutions to Markovian superquadratic BSDEs with an unbounded terminal condition}

\author{
Federica Masiero\\
Dipartmento di Matematica e Applicazioni, Università di Milano Bicocca\\
via Cozzi 53, 20125 Milano, Italy\\
e-mail: federica.masiero@unimib.it\\
\\
Adrien Richou\\
Univ. Bordeaux, IMB, UMR 5251, F-33400 Talence, France.\\
CNRS, IMB, UMR 5251, F-33400 Talence, France.\\
INRIA, \'Equipe ALEA, F-33400 Talence, France.\\ 
e-mail: adrien.richou@math.u-bordeaux1.fr}

\selectlanguage{english}

\maketitle

\begin{abstract} 
In \cite{Richou-12}, the author proved the existence and the uniqueness of solutions to Markovian superquadratic BSDEs with an unbounded terminal condition when the generator and the terminal condition are locally Lipschitz. In this paper, we prove that the existence result remains true for these BSDEs when the regularity assumptions on the terminal condition is weakened.
\end{abstract}

% \selectlanguage{francais}
% 
% \begin{abstract} 
% 
% \end{abstract}
\selectlanguage{english}

%\paragraph{Key words and phrases.} 

%\paragraph{AMS subject classifications.} 60H10, 60H30, 93E20.
%\paragraph{AMS subject classifications.} 60H10.

\section{Introduction}

Since the early nineties and the work of Pardoux and Peng \cite{Pardoux-Peng-90}, there has been an increasing interest for backward stochastic differential equations (BSDEs for short) because of the wide range of applications. A particular class of BSDE is studied since few years: BSDEs with generators of quadratic growth with respect to the variable $z$ (quadratic BSDEs for short). See e.g. \cite{Kobylanski-00,Briand-Hu-06,Delbaen-Hu-Richou-09} for existence and uniqueness results and \cite{Rouge-ElKaroui-00,Hu-Imkeller-Muller-05,Mania-Schweizer-05} for applications.

Naturally, we could also wonder what happens when the generator has a superquadratic growth with respect to the variable $z$. Up to our knowledge the case of superquadratic BSDEs was firstly investigated in the recent paper \cite{Delbaen-Hu-Bao-09}. In this article, the authors consider superquadratic BSDEs when the terminal condition is bounded and the generator is convex in $z$. Firstly, they show that in a general way the problem is ill-posed: given a superquadratic generator, there exists a bounded terminal condition such that the associated BSDE does not admit any bounded solution and, on the other hand, if the BSDE admits a bounded solution, there exist infinitely many bounded solutions for this BSDE. In the same paper, the authors also show that the problem becomes well-posed in a Markovian framework: when the terminal condition and the generator are deterministic functions of a forward SDE, we have an existence result. More precisely, let us consider $(X,Y,Z)$ the solution to the (decoupled) forward 
backward 
system
\begin{eqnarray*}
 X_t &=& x +\int_0^t b(s,X_s)ds+\int_0^t \sigma(s) dW_s,\\
 Y_t &=& g(X_T) +\int_t^T f(s,X_s,Y_s,Z_s)ds-\int_t^T Z_s dW_s,
 \end{eqnarray*}
with growth assumptions
\begin{eqnarray*}
 \abs{f(t,x,y,z)} &\leqslant& C(1+\abs{x}^{p_f}+\abs{y}+\abs{z}^{l+1}), \quad l>1,\\
 \abs{g(x)} &\leqslant& C(1+\abs{x}^{p_g}).
\end{eqnarray*}
In \cite{Delbaen-Hu-Bao-09}, the authors obtain an existence result by assuming that $p_g=p_f=0$, $f$ is a convex function that depends only on $z$ and $g$ is a lower (or upper) semi-continuous function. As in the quadratic case it is possible to show that the boundedness of the terminal condition is a too strong assumption: in \cite{Richou-12}, the author shows an existence and uniqueness result by assuming that $p_g \leqslant 1+1/l$, $p_f \leqslant 1+1/l$, $f$ and $g$ are locally Lipschitz functions with respect to $x$ and $z$. When we consider this result, two questions 
arise:
\begin{itemize}
 \item Could we have an existence result when $p_g$ or $p_f$ is greater than $1+1/l$ ?
 \item Could we have an existence result when $f$ or $g$ is less smooth with respect to $x$ or $z$, that is to say, is it possible to have assumptions on the growth of $g$ and $f$ but not on the growth of their derivatives with respect to $x$ and $z$ ?
\end{itemize}
For the first question, the answer is clearly ``no'' in the quadratic case: see e.g. \cite{Delbaen-Hu-Richou-09}. In the superquadratic case, the authors of \cite{Gladkov-Guedda-Kersner-08} have obtained the same limitation on the growth of the initial condition for the so-called generalized deterministic KPZ equation $u_t=u_{xx}+\lambda\abs{u_x}^q$ and they show that this boundary is sharp for power-type initial conditions. So, it seems that the answer of the first question is also ``no'' in the superquadratic case.

For the second question, the answer is clearly ``yes'' in the quadratic case. Indeed, a smoothness assumption on $f$ is required for uniqueness results (see e.g. \cite{Briand-Hu-08,Delbaen-Hu-Richou-09}) but not for existence results (see e.g. \cite{Briand-Hu-08,Barrieu-ElKaroui-11}). In the superquadratic case, the authors of \cite{Delbaen-Hu-Bao-09} show an existence result when $g$ is only lower (or upper) semi-continuous but also bounded. Nevertheless $f(z)$ is assumed to be convex, that implies that it is a locally Lipschitz function. The aim of this note is to mix results of articles \cite{Delbaen-Hu-Bao-09,Richou-12} to obtain an existence result when the terminal condition is only lower (or upper) semi-continuous and unbounded. Let us remark that we answer only partially to the second question because we do not relax smoothness assumptions on $f$.

For completeness, in the recent paper \cite{Cheridito-Stadje-12}, Cheridito and Stadje show an existence and uniqueness result for superquadratic BSDEs in a Lipschitz or bounded ``path-dependent'' framework: the terminal condition and the generator are Lipschitz or bounded functions of Brownian motion paths. To the best of our knowledge, \cite{Delbaen-Hu-Bao-09,Richou-12,Cheridito-Stadje-12} are the only papers that deal with superquadratic BSDEs.
 
The paper is organized as follows. In section 2 we obtain some general a priori estimates on $Y$ and $Z$ for Markovian superquadratic BSDEs whereas section 3 is devoted to the existence result described before.

\paragraph{Notations}
Throughout this paper, $(W_t)_{t \geqslant 0}$ will denote a $d$-dimensional Brownian motion, defined on a probability space $(\Omega,\mathcal{F}, \P)$. For $t \geqslant 0$, let $\mathcal{F}_t$ denote the $\sigma$-algebra $\sigma(W_s; 0\leqslant s\leqslant t)$, augmented with the $\P$-null sets of $\mathcal{F}$. The Euclidean norm on $\R^d$ will be denoted by $|.|$. The operator norm induced by $|.|$ on the space of linear operators is also denoted by $|.|$. The notation $\mathbb{E}_t$ stands for the conditional expectation given $\mathcal{F}_t$. For $p \geqslant 2$, $m \in \N$, we denote further
\begin{itemize}
 \item $\mathcal{S}^p$ the space of real-valued, adapted and càdlàg processes $(Y_t)_{t \in [0,T]}$ normed by $\norm{Y}_{\mathcal{S}^p}=\E [(\sup_{t \in [0,T]} \abs{Y_t})^p]^{1/p}$;
 \item $\mathcal{M}^p(\R^m)$, or $\mathcal{M}^p$, the space of all progressively measurable processes $(Z_t)_{t \in [0,T]}$ with values in $\R^m$ normed by $\norm{Z}_{\mathcal{M}^p}=\E[(\int_0^T \abs{Z_s}^2ds)^{p/2}]^{1/p}$.
\end{itemize}
 In the following, we keep the same notation $C$ for all finite, nonnegative constants that appear in our computations.%: they may depend on known parameters deriving from assumptions and on $T$, but not on any of the approximation and discretization parameters.

In this paper we consider $X$ the solution to the SDE
\begin{equation}
\label{EDS}
 X_t=x+\int_0^t b(s,X_s)ds+\int_0^t \sigma(s) dW_s,
\end{equation}
and $(Y,Z) \in \mathcal{S}^2\times \mathcal{M}^2$ the solution to the Markovian BSDE
\begin{equation}
\label{EDSR}
 Y_t=g(X_T)+\int_t^T f(s,X_s,Y_s,Z_s)ds-\int_t^T Z_sdW_s.
\end{equation}
By a solution to the BSDE (\ref{EDSR}) we mean a pair $(Y_t,Z_t)_{t \in [0,T]}$ of predictable processes with values in $\R \times \R^{1 \times d}$ such that $\P$-a.s., $t \mapsto Y_t$ is continuous, $t \mapsto Z_t$ belongs to $L^2([0,T])$, $t \mapsto f(t,X_t,Y_t,Z_t)$ belongs to $L^1([0,T])$ and $\mathbb{P}-a.s.$ the equation (\ref{EDSR}) is verified.

\section{Some a priori estimates on $Y$ and $Z$}
For the SDE (\ref{EDS}) we use standard assumption.
\paragraph{Assumption (F.1).}
Let $b : [0,T] \times \mathbb{R}^d \rightarrow \mathbb{R}^d$ and $\sigma : [0,T] \rightarrow \mathbb{R}^{d \times d}$ be continuous  functions and let us assume that there exists  $K_b \geqslant 0$ such that:
\begin{enumerate}[(a)]
 \item $\forall t \in [0,T]$, $\abs{b(t,0)} \leqslant C$,
 \item $\forall t \in [0,T]$, $\forall (x,x') \in \mathbb{R}^d \times \mathbb{R}^d$, $\abs{b(t,x)-b(t,x')} \leqslant K_b \abs{x-x'}.$
\end{enumerate}
Let us now consider the following assumptions on the generator and on the terminal condition of the BSDE (\ref{EDSR}).
\paragraph{Assumption (B.1).}
Let $f: [0,T] \times \mathbb{R}^d \times \mathbb{R} \times \mathbb{R}^{1\times d} \rightarrow \mathbb{R}$ be a continuous function and let us assume that there exist five constants, $l > 1$, $0 \leqslant r_f < \frac{1}{l}$, $\beta \geqslant 0$, $\gamma \geqslant 0$ and $\delta \geqslant 0$ such that:
\begin{enumerate}[(a)]
\item for each $(t,x,y,y',z) \in [0,T] \times \mathbb{R}^d \times \mathbb{R} \times \mathbb{R} \times \mathbb{R}^{1\times d}$,
$$ \abs{f(t,x,y,z)-f(t,x,y',z)} \leqslant \delta \abs{y-y'};$$
\item for each $(t,x,y,z,z') \in [0,T] \times \R^d \times \R \times \R^{1\times d} \times \R^{1\times d}$,
$$\abs{ f(t,x,y,z)-f(t,x,y,z')} \leqslant \left(C+\frac{\gamma}{2}(\abs{z}^l+\abs{z'}^l)\right)\abs{z-z'};$$
\item for each $(t,x,x',y,z) \in [0,T] \times \R^d \times \R^d \times \R \times \R^{1\times d}$,
$$\abs{ f(t,x,y,z)-f(t,x',y,z)} \leqslant \left(C+\frac{\beta}{2}(\abs{x}^{r_f}+\abs{x'}^{r_f})\right)\abs{x-x'}.$$
\end{enumerate}
\paragraph{Assumption (TC.1).}
Let $g:\mathbb{R}^d \rightarrow  \mathbb{R}$ be a continuous function and let us assume that there exist $0 \leqslant r_g < \frac{1}{l}$ and $\alpha \geqslant 0$ such that: for each $(t,x,x',y,z) \in [0,T] \times \R^d \times \R^d \times \R \times \R^{1\times d}$,
$$\abs{g(x)-g(x')} \leqslant \left(C+\frac{\alpha}{2}(\abs{x}^{r_g}+\abs{x'}^{r_g})\right)\abs{x-x'}.$$

We also use more general growth assumptions that are more natural for existence results.  
\paragraph{Assumptions (B.2).}
Let $f: [0,T] \times \mathbb{R}^d \times \mathbb{R} \times \mathbb{R}^{1\times d} \rightarrow \mathbb{R}$ be a continuous function and let us assume that there exist constants, $l > 1$, $0 \leqslant r_f < \frac{1}{l}$, $\bar{\beta} \geqslant 0$, $\bar{\gamma} \geqslant 0$, $\bar{\delta} \geqslant 0$, $0\leqslant \eta < l+1$, $\varepsilon>0$ such that: one of these inequalities holds, for all $(t,x,y,z) \in [0,T] \times \mathbb{R}^d \times \mathbb{R} \times \mathbb{R}^{1\times d}$, 
 \begin{enumerate}[(a)]
  \item $\abs{f(t,x,y,z)} \leqslant C+\bar{\beta}\abs{x}^{r_f+1}+\bar{\delta}\abs{y}+\bar{\gamma} \abs{z}^{l+1}$,
  \item $-C-\bar{\beta} \abs{x}^{r_f+1}-\bar{\delta}\abs{y}-\bar{\gamma} \abs{z}^{\eta} \leqslant f(t,x,y,z) \leqslant C+\bar{\beta}\abs{x}^{r_f+1}+\bar{\delta}\abs{y}+\bar{\gamma} \abs{z}^{l+1}$,
  %\item $-C-\bar{\beta} \abs{x}^{r+1}-\bar{\delta}\abs{y}-\bar{\gamma} \abs{z}^2 \leqslant f(t,x,y,z) \leqslant C+\bar{\beta}\abs{x}^{r+1}+\bar{\delta}\abs{y}+\bar{\gamma} \abs{z}^{l+1}$,
  \item $-C-\bar{\beta} \abs{x}^{r_f+1}-\bar{\delta}\abs{y}+\varepsilon \abs{z}^{l+1} \leqslant f(t,x,y,z) \leqslant C+\bar{\beta}\abs{x}^{r_f+1}+\bar{\delta}\abs{y}+\bar{\gamma} \abs{z}^{l+1}$.
 \end{enumerate}

\paragraph{Assumption (TC.2).}
Let $g:\mathbb{R}^d \rightarrow  \mathbb{R}$ be a lower semi-continuous function and let us assume that there exist $0 \leqslant p_g < 1+1/l$ and $\bar{\alpha} \geqslant 0$ such that: for each $x \in \mathbb{R}^d$,
$$\abs{g(x)} \leqslant C+\bar{\alpha} \abs{x}^{p_g}.$$

\begin{rem}
The following relations hold true:
\begin{itemize}
 \item (B.2)(c) $\Rightarrow$ (B.2)(b)  $\Rightarrow$ (B.2)(a).
 \item (B.1)  $\Rightarrow$ (B.2)(a).
 \item (TC.1)  $\Rightarrow$ (TC.2) with $p_g=r_g+1$.
 \item We only consider superquadratic BSDEs, so $l > 1$. $l=1$ corresponds to the quadratic case.
\end{itemize}
\end{rem}

Firstly, let us recall the existence and uniqueness result shown in \cite{Richou-12}.

\begin{prop}
\label{existence unicite localement lipschitz}
We assume that (F.1), (B.1) and (TC.1) hold. There exists a solution $(Y,Z)$ of the Markovian BSDE (\ref{EDSR}) in $\mathcal{S}^2\times \mathcal{M}^2$ such that,
\begin{equation} 
\label{estimee croissance Z}
\abs{Z_t} \leqslant A+B(\abs{X_t}^{r_g}+(T-t)\abs{X_t}^{r_f}), \quad \forall t \in [0,T].
\end{equation}
Moreover, this solution is unique amongst solutions $(Y,Z)$ such that
\begin{itemize}
 \item $Y \in \mathcal{S}^2$,
 \item there exists $\eta >0$ such that
$$\mathbb{E} \left[e^{(\frac{1}{2}+\eta)\frac{\gamma^2}{4}\int_0^T \abs{Z_s}^{2l} ds}\right] < +\infty.$$
\end{itemize}
\end{prop}
\begin{rem}
 To be precise, in the Proposition 2.2 of the article \cite{Richou-12} the author shows the estimate
$$\abs{Z_t} \leqslant A+B\abs{X_t}^{r_g \vee r_f}, \quad \forall t \in [0,T],$$
but it is rather easy to do the proof again to show the estimate (\ref{estimee croissance Z}) given in Proposition \ref{existence unicite localement lipschitz}.
\end{rem}

Such a result allows us to obtain a comparison result.
\begin{prop}
\label{comparison result}
 We assume that (F.1) holds. Let $f_1$, $f_2$ two generators and $g_1$, $g_2$ two terminal conditions such that (B.1) and (TC.1) hold. Let $(Y^1,Z^1)$ and $(Y^2,Z^2)$ be the associated solutions given by Proposition \ref{existence unicite localement lipschitz}. We assume that $g_1 \leqslant g_2$ and $f_1 \leqslant f_2$. Then we have that $Y^1 \leqslant Y^2$ almost surely.
\end{prop}

\paragraph*{Proof of the proposition}
The proof is the same than the classical one that can be found in \cite{ElKaroui-Peng-Quenez-97} for example. Let us set $\delta Y:=Y^1-Y^2$ and $\delta Z:=Z^1-Z^2$. The usual linearization trick gives us
$$\delta Y_t=g_1(X_T)-g_2(X_T)+\int_t^T f_1(s,X_s,Y^1_s,Z^1_s)-f_2(s,X_s,Y^1_s,Z^1_s)+\delta Y_s U_s +\delta Z_s V_s ds -\int_t^T \delta Z_s dW_s,$$
with $\abs{U_s}\leqslant \delta$ and 
$$\abs{V_s} \leqslant C+\frac{\gamma}{2}\left( \abs{Z^1_s}^l +\abs{Z^2_s}^l\right) \leqslant C(1+\abs{X_s}^{(r_g\vee r_f)l}).$$
Since $(r_g\vee r_f)l<1$, Novikov's condition is fulfilled and we are allowed to apply Girsanov's transformation:
\begin{eqnarray*}
 \delta Y_t &=& \E^{\Q}_t \left[ e^{\int_t^T U_udu}(g_1(X_T)-g_2(X_T))+\int_t^T e^{\int_t^s U_udu}(f_1(s,X_s,Y^1_s,Z^1_s)-f_2(s,X_s,Y^1_s,Z^1_s))ds \right]\\
&\leqslant& 0,
\end{eqnarray*}
with
$$\frac{d\Q}{d\P}=\exp \left( \int_0^T V_s dW_s-\frac{1}{2} \int_0^T \abs{V_s}^2ds \right).$$
\cqfd

Now we are ready to prove estimates on $Y$ and $Z$.

\begin{prop}
\label{estimation Y}
 Let us assume that (F.1), (B.1), (B.2), (TC.1) and (TC.2) hold. Let $(Y,Z)$ be the solution of the BSDE (\ref{EDSR}) given by Proposition \ref{existence unicite localement lipschitz}. Then we have, for all $t \in [0,T]$,
$$\abs{Y_t} \leqslant C(1+\abs{X_t}^{p_g}+(T-t)\abs{X_t}^{r_f+1})$$
with a constant $C$ that depends on constants that appear in assumptions (F.1), (B.2) and (TC.2) but not in assumptions (B.1) and (TC.1).
\end{prop}
\paragraph*{Proof of the proposition}
Let us consider the terminal condition
$$\bar{g}(x)=C+\bar{\alpha} (\abs{x}+1)^{p_g},$$
and the generator
$$\bar{f}(t,x,y,z) =C+\bar{\beta}\abs{x}^{r_f+1}+\bar{\delta}\abs{y}+\bar{\gamma} \abs{z}^{l+1},$$
with $C$ such that $g \leqslant \bar{g}$ and $f \leqslant \bar{f}$. (B.1) holds for $\bar{f}$ and (TC.1) holds for $\bar{g}$, so, according to Proposition \ref{existence unicite localement lipschitz}, there exists a unique solution $(\bar{Y},\bar{Z})$ to the BSDE
$$ \bar{Y}_t=\bar{g}(X_T)+\int_t^T \bar{f}(s,X_s,\bar{Y}_s,\bar{Z}_s)ds-\int_t^T \bar{Z}_sdW_s.$$
Thanks to Proposition \ref{comparison result}, we know that
$$Y \leqslant \bar{Y}, \quad \textrm{and} \quad \bar{Y}\geqslant 0.$$
Moreover, since $\abs{\bar{Z}_s} \leqslant C(1+\abs{X_s}^{(p_g-1)\vee r_f})$, $(p_g-1)l<1$ and $r_f l<1$, we have
\begin{eqnarray*}
\bar{Y}_t &\leqslant&\E_t \left[e^{\bar{\delta}(T-t)} (C+\bar{\alpha}(\abs{X_t}+1)^{p_g})+\int_t^T e^{\bar{\delta}(s-t)}(C+\bar{\beta}\abs{X_s}^{r_f+1}+\bar{\gamma}\abs{\bar{Z}_s}^{l+1})ds\right]\\
&\leqslant & C\left(1+\E_t\left[\sup_{t \leqslant s\leqslant T}\abs{X_s}^{p_g}\right]+(T-t)\E_t\left[\sup_{t \leqslant s\leqslant T}\abs{X_s}^{r_f+1}\right]\right).
\end{eqnarray*}
Let us remark that the constant $C$ in the a priori estimate for $\bar{Z}$ depends on constants that appear in assumptions (F.1), (B.2) and (TC.2) but not in assumptions (B.1) and (TC.1). Thanks to classical estimates on SDEs we have, for all $p \geqslant 1$,
$$\E_t\left[\sup_{t \leqslant s\leqslant T}\abs{X_s}^{p}\right] \leqslant C(1+\abs{X_t}^{p}),$$
so we obtain
$$Y_t \leqslant \bar{Y}_t \leqslant C(1+\abs{X_t}^{p_g}+(T-t)\abs{X_t}^{r_f+1}).$$
By the same type of argument we easily show that
$$ -C(1+\abs{X_t}^{p_g}+(T-t)\abs{X_t}^{r_f+1})\leqslant Y_t,$$
and this concludes the proof.
\cqfd

\begin{prop}
\label{estimation Z}
 Let us assume that (F.1), (B.1), (B.2)(c), (TC.1) and (TC.2) hold. Let $(Y,Z)$ be the solution of the BSDE (\ref{EDSR}) given by Proposition \ref{existence unicite localement lipschitz}. Then, for all $t \in [0,T]$, we have 
 $$\E_t\left[\int_t^T\abs{Z_s}^{l+1}ds\right] \leqslant C(1+\abs{X_t}^{p_g}+(T-t)\abs{X_t}^{r_f+1}),$$
with a constant $C$ that depends on constants that appear in assumptions (F.1), (B.2)(c) and (TC.2) but not in assumptions (B.1) and (TC.1).
\end{prop}

\paragraph*{Proof of the proposition}
To show the proposition we just have to write
\begin{eqnarray*}
 \E_t\left[\int_t^T\abs{Z_s}^{l+1}ds\right] &\leqslant& \frac{1}{\varepsilon}\left(\E_t\left[\int_t^T f(s,X_s,Y_s,Z_s)ds +\int_t^T \left( C+\bar{\beta} \abs{X_s}^{r_f+1}+\bar{\delta}\abs{Y_s}\right)ds\right]\right)\\
&\leqslant& \frac{1}{\varepsilon}\left(\E_t\left[Y_t-g(X_T) +\int_t^T C+\bar{\beta} \abs{X_s}^{r_f+1}+\bar{\delta}\abs{Y_s}ds\right]\right)\\
&\leqslant& C(1+(T-t)\abs{X_t}^{r_f+1}+\abs{X_t}^{p_g})
\end{eqnarray*}
thanks to Proposition \ref{estimation Y}. \cqfd

\begin{rem}
\label{inversion hypotheses croissance f}
Proposition \ref{estimation Z} stays true if we replace assumption (B.2)(c) by
$$-C-\bar{\beta} \abs{x}^{r_f+1}-\bar{\delta}\abs{y}- \bar{\gamma} \abs{z}^{l+1} \leqslant f(t,x,y,z) \leqslant C+\bar{\beta}\abs{x}^{r_f+1}+\bar{\delta}\abs{y}-\varepsilon\abs{z}^{l+1}.$$
\end{rem}

\begin{rem}
 In Propositions \ref{estimation Y} and \ref{estimation Z} we insist on the fact that $C$ does not depend on constants that appear in assumptions (B.1) and (TC.1) when the local Lipschitzianity of the coefficients is stated. Thanks to this property, we can use these a priori estimates on $Y$ and $Z$ in the following section where we obtain an existence result when the terminal condition is not locally Lipschitz.
 
\end{rem}

\section{An existence result}
Let us now introduce new assumptions.
\paragraph{Assumption (F.2).}
$b$ is differentiable with respect to $x$ and $\sigma$ is differentiable with respect to $t$. There exists $\lambda \in \R^+$ such that
$\forall \eta \in \mathbb{R}^d$ 
\begin{equation*}
\label{hypothese sur nabla b}
\abs{\tr{\eta}\sigma(s)[\tr{\sigma(s)}\tr{\nabla b(s,x)}-\tr{\sigma'(s)}]\eta} \leqslant \lambda\abs{\tr{\eta}\sigma(s)}^2, \quad \forall (s,x) \in [0,T] \times \R^d.
\end{equation*}
\begin{rem}
It is shown in part 5.5.1 of \cite{Richou-10} that if $\sigma$ does not depend on time, assumption (F.2) is equivalent to this kind of commutativity assumption:
\begin{itemize}
 \item there exist $A : [0,T] \times \R^d \rightarrow \R^{d \times d}$ and $B : [0,T] \rightarrow \R^{d \times d}$ such that $A$ is differentiable with respect to $x$, $\nabla_x A$ is bounded and $\forall x \in \R^d$, $\forall s \in [0,T]$, $b(s,x)\sigma=\sigma A(s,x)+B(s).$ 
\end{itemize}
It is also noticed in \cite{Richou-10} that this assumption allows us to reduce assumption on the regularity of $b$ by a standard smooth approximation of $A$. 
\end{rem}

\paragraph{Assumption (B.3).}
$f$ is differentiable with respect to $z$ and for all $(t,x,y,z) \in [0,T] \times \mathbb{R}^d \times \mathbb{R} \times \mathbb{R}^{1\times d}$,
$$f(t,x,y,z)-\langle \nabla_zf(t,x,y,z),z\rangle \leqslant C-\varepsilon \abs{z}^{l+1}.$$

\begin{rem}
Let us give some substantial examples of functions such that (B.3) holds. If we assume that $f(t,x,y,z):=f_1(t,x,y,z)+f_2(t,x,y,z)$ with $f_1$ a differentiable function with respect to $z$ such that, $\exists p \in[0,l[$,  $\forall (t,x,y,z) \in [0,T] \times \mathbb{R}^d \times \mathbb{R} \times \mathbb{R}^{1\times d}$,
$$\abs{\nabla_z f_1(t,x,y,z)} \leqslant (1+\abs{z}^{p}),$$
and $f_2$ is a twice differentiable function with respect to $z$ such that,  $\forall (t,x,y,z) \in [0,T] \times \mathbb{R}^d \times \mathbb{R} \times \mathbb{R}^{1\times d}$, $\forall u \in \R^d$,
$$^tu\nabla^2_{zz}f_2(t,x,y,z)u \geqslant (-C+\varepsilon\abs{z}^{l-1})\abs{u}^2,$$
then we easily see that 
$$f_1(t,x,y,z)-\langle \nabla_zf_1(t,x,y,z),z\rangle \leqslant C+C \abs{z}^{p+1},$$
and a direct application of Taylor expansion with integral form gives us
$$f_2(t,x,y,z)-\langle \nabla_zf_2(t,x,y,z),z\rangle \leqslant C-C' \abs{z}^{l+1},$$
so (B.3) holds. For example, (B.3) holds for the function $z \mapsto C\abs{z}^{l+1}+h(\abs{z}^{l+1-\eta})$ with $C>0$, $0<\eta\leqslant l+1$ and $h$ a differentiable function with a bounded derivative.
\end{rem}

\begin{prop}
\label{prop estimation temporelle Z}
 Let us assume that (F.1), (F.2), (B.1), (B.3), (TC.1) and (TC.2) hold. Let $(Y,Z)$ be the solution of the BSDE (\ref{EDSR}) given by Proposition \ref{existence unicite localement lipschitz}. If we assume that $0 \leqslant p_gl <1$, then we have, for all $t \in [0,T[$,
 $$\abs{Z_t} \leqslant \frac{C(1+\abs{X_t}^{p_g/(l+1)})}{(T-t)^{1/(l+1)}}+C\abs{X_t}^{\frac{r_f+1}{l+1} }.$$
The constant $C$ depends on constants that appear in assumptions (F.1), (F.2), (B.1), (B.3) and (TC.2) but not in assumption (TC.1).
\end{prop}

\paragraph*{Proof of the proposition}
Firstly we approximate our Markovian BSDE by another one. Let $(Y^M,Z^M)$ the solution of the BSDE
\begin{equation}
\label{EDSR approchee}
Y^M_t = g_M(X_T)+\int_t^T f_M(s,X_s,Y^M_s,Z_s^M)ds-\int_t^T Z_s^M dW_s,
\end{equation}
with $g_M=g \circ \rho_M$ and $f_M=f(.,\rho_M(.),.,.)$ where $\rho_M$ is a smooth modification of the projection on the centered Euclidean ball of radius $M$ such that $\abs{\rho_M}\leqslant M$, $\abs{\nabla \rho_M} \leqslant 1$ and $\rho_M(x)=x$ when $\abs{x}\leqslant M-1$. It is now easy to see that $g_M$ and $f_M$ are Lipschitz functions with respect to $x$. Proposition 2.3 in \cite{Richou-12} gives us that $Z^M$ is bounded by a constant $C_0$ that depends on $M$. So, $f_M$ is a Lipschitz function with respect to $z$ and BSDE (\ref{EDSR approchee}) is a classical Lipschitz BSDE. Now we use the following Lemma that will be shown afterwards.
\begin{lem}
\label{lemme recurrence}
Let us assume that (F.1), (F.2), (B.1), (B.3), (TC.1) and (TC.2) hold. We also assume that $0 \leqslant p_gl <1$. Then we have, for all $t \in [0,T[$,
 $$\abs{Z_t^M} \leqslant \frac{A_n+B_n\abs{X_t}^{p_g/(l+1)}}{(T-t)^{1/(l+1)}}+D_n\abs{X_t}^{\frac{r_f+1}{l+1} },$$
with $(A_n,B_n,D_n)_{n \in \N}$ defined by recursion: $B_0=0$, $D_0=0$, $A_0=C_0T^{1/(l+1)}$,
$$A_{n+1}=C(1+A_n^{al}+B_n^{alp}+D_n^{al\bar{p}}), \quad B_{n+1}=C, \quad D_{n+1}=C,$$
where $a:=(p_g\vee (r_f+1))/(l+1)$, $p>1$, $\bar{p}>1$ and $C$ is a constant that does not depend on $M$ and constants in assumption (TC.1).
\end{lem}
Since $al<1$, the recursion function that define the sequence $(A_n)_{n \geqslant 0}$ is a contractor function, so $A_n \rightarrow A_{\infty}$ when $ n \rightarrow +\infty$, with $A_{\infty}$ that does not depend on $M$ and constants in assumption (TC.1). Finally, we have, for all $t \in [0,T[$,
 $$\abs{Z_t^M} \leqslant \frac{C(1+\abs{X_t}^{p_g/(l+1)})}{(T-t)^{1/(l+1)}}+C\abs{X_t}^{\frac{r_f+1}{l+1}}.$$
The constant $C$ depends on constants that appear in assumptions (F.1), (F.2), (B.1), (B.3) and (TC.2) but not in assumption (TC.1). Moreover $C$ does not depends on $M$. Now, we want to come back to the initial BSDE (\ref{EDSR}). It is already shown in the proof of Proposition 2.2 of the article \cite{Richou-12} that $(Y^n,Z^n) \rightarrow (Y,Z)$ in $\mathcal{S}^2 \times \mathcal{M}^2$. So our estimate on $Z^M$ stays true for a version of $Z$.
\cqfd

\paragraph*{Proof of Lemma \ref{lemme recurrence}}
Let us prove the result by recursion. For $n=0$ we have already shown the result. Let us assume that the result is true for some $n \in \N$ and let us show that it stays true for $n+1$. In a first time we suppose that $f$ and $g$ are differentiable with respect to $x$ and $y$. Then $(Y^M,Z^M)$ is differentiable with respect to $x$ and $(\nabla Y^M,\nabla Z^M)$ is the solution of the BSDE
\begin{eqnarray*}
 \nabla Y_t^M &=& \nabla g_M(X_T)\nabla X_T - \int_t^T \nabla Z_s^M dW_s\\
& & +\int_t^T \nabla_x f_M(s,X_s,Y_s^M,Z_s^M) \nabla X_s + \nabla_y f_M(s,X_s,Y_s^M,Z_s^M) \nabla Y_s^M + \nabla_z f_M(s,X_s,Y_s^M,Z_s^M) \nabla Z_s^M ds,
\end{eqnarray*} 
and a version of $Z^M$ is given by  $(\nabla Y_t^M(\nabla X_t)^{-1} \sigma(t))_{t \in[0,T]}$. Let us introduce some notations: we set
\begin{eqnarray*}
 d\tilde{W}_t &:=& dW_t-\nabla_z f_M(t,X_t,Y_t^M,Z_t^M)dt,\\
\alpha_t &:= & \int_0^t e^{\int_0^s\nabla_y f_M(u,X_u,Y_u^M,Z_u^M)du}\nabla_x f_M(s,X_s,Y_s^M,Z_s^M)\nabla X_s ds (\nabla X_t)^{-1}\sigma(t),\\
\tilde{Z}_t^M &:=&  e^{\int_0^t\nabla_y f_M(s,X_s,Y_s^M,Z_s^M)ds} Z_t^M +\alpha_t.
\end{eqnarray*}
By applying Girsanov's theorem we know that there exists a probability $\Q^M$ under which $\tilde{W}$ is a Brownian motion with
$$\frac{d\Q^M}{d\P} = \exp\left( \int_0^T \nabla_z f_M(t,X_t,Y_t^M,Z_t^M)dW_t-\frac{1}{2}\int_0^T \abs{\nabla_z f_M(t,X_t,Y_t^M,Z_t^M)}^2dt\right).$$
Then, exactly as in the proof of Theorem 3.3 in \cite{Richou-11}, we can show the following lemma. 
\begin{lem}
\label{Qm sous martingale}
$\abs{e^{\lambda t} \tilde{Z}_t^M}^2$ is a $\mathbb{Q}^M$-submartingale. 
\end{lem}
For the reader's convenience, we recall this proof in the appendix. It results that $\abs{e^{\lambda t} \tilde{Z}_t^M}^{l+1}$ is also a $\mathbb{Q}^M$-submartingale and we have:
\begin{eqnarray*}
 \mathbb{E}^{\mathbb{Q}^M}_t \left[\int_t^T e^{2\lambda s}\abs{\tilde{Z}^M_s}^{l+1}ds \right] & \geqslant & e^{2\lambda t}\abs{\tilde{Z}_t^M}^{l+1}(T-t)\\
&\geqslant & e^{2\lambda t}\abs{e^{\int_0^t\nabla_y f_M(s,X_s,Y_s^M,Z_s^M)ds}Z_t^M+\alpha_t}^{l+1}(T-t) ,
\end{eqnarray*}
which implies
\begin{eqnarray}
\nonumber
  \abs{Z_t^M}^{l+1}(T-t) & \leqslant & C\left(e^{2\lambda t}\abs{e^{\int_0^t\nabla_y f_M(s,X_s,Y_s^M,Z_s^M)ds}Z_t^M+\alpha_t}^{l+1}+\abs{\alpha_t}^{l+1}\right)(T-t)\\
\nonumber &\leqslant & C\left( \mathbb{E}^{\mathbb{Q}^M}_t \left[\int_t^T e^{2\lambda s}\abs{\tilde{Z}_s^M}^{l+1}ds\right] +(T-t)\left(1+\abs{X_t}^{(l+1)r_f}\right)\right)\\
\nonumber &\leqslant & C\left( 1+\mathbb{E}^{\mathbb{Q}^M}_t \left[\int_t^T \abs{Z_s^M}^{l+1}ds\right] +\mathbb{E}^{\mathbb{Q}^M}_t \left[\int_t^T \abs{X_s}^{(l+1)r_f}ds\right]+(T-t)\abs{X_t}^{(l+1)r_f}\right).\\
\label{inegalite1 Z} & &
\end{eqnarray}
Let us recall that $(Y^M,Z^M)$ is solution of BSDE
$$Y_t^M=g_M(X_T)+\int_t^T \tilde{f}_M(s,X_s,Y_s^M,Z_s^M)ds-\int_t^T Z_s^Md\tilde{W}_s,$$
with
$$\tilde{f}_M(s,x,y,z):=f_M(s,x,y,z)-\langle z,\nabla_zf_M(s,x,y,z)\rangle.$$
Since assumption (B.3) holds for $f$, assumption (B.2)(c) holds for $-\tilde{f}_M$ with constants that do not depend on $M$. Then we can mimic the proof of Proposition \ref{estimation Z} (see also Remark \ref{inversion hypotheses croissance f}) to show that
\begin{equation}
\label{inegalite2 Z}
\mathbb{E}^{\mathbb{Q}^M}_t \left[\int_t^T \abs{Z_s^M}^{l+1}ds\right] \leqslant C\left(1+\mathbb{E}^{\mathbb{Q}^M}_t \left[\abs{X_T}^{p_g}\right]+\int_t^T \left(\mathbb{E}^{\mathbb{Q}^M}_t \left[\abs{X_s}^{p_g}\right]+\mathbb{E}^{\mathbb{Q}^M}_t \left[\abs{X_s}^{r_f+1}\right]\right)ds\right),
\end{equation}
with a constant $C$ that does not depend on $M$ and constants that appear in assumption (TC.1). Then, by putting (\ref{inegalite2 Z}) in (\ref{inegalite1 Z}), we see that we just have to obtain an a priori estimate for $\mathbb{E}^{\mathbb{Q}^M}_t \left[\abs{X_s}^{c}\right]$ with $c \in \R^{+*}$. We have 
\begin{eqnarray*}
 \abs{X_s} &=& \abs{X_t+\int_t^s b(u,X_u)du+\int_t^s \sigma(u) d\tilde{W}_u+\int_t^s \sigma(u) \nabla_z f_M(u,X_u,Y_u^M,Z_u^M)du}\\
&\leqslant& \abs{X_t}+C+C\int_t^s\abs{X_u}du+\abs{\int_t^s \sigma(u) d\tilde{W}_u}+C\int_t^s \abs{Z_u^M}^ldu,
\end{eqnarray*}
with $C$ that does not depend on $M$. Now we use the recursion assumption to obtain
\begin{eqnarray*}
 \int_t^s \abs{Z_u^M}^ldu & \leqslant& C\int_t^s \left(\frac{A_n^l}{(T-u)^{l/(l+1)}} +\frac{B_n^l}{(T-u)^{l/(l+1)}}\abs{X_u}^{lp_g/(l+1)}+D_n^l\abs{X_u}^{(r_f+1)l/(l+1)}\right)du.
\end{eqnarray*}
Obviously we have $\int_t^T \frac{A_n^l}{(T-u)^{l/(l+1)}}du \leqslant CA_n^l$. For the other terms we use Young inequality: Since $lp_g/(l+1) <1$ and $(r_f+1)l/(l+1)<1$ , we have
\begin{eqnarray*}
 \int_t^s \abs{Z_u^M}^ldu & \leqslant& CA_n^l+C\int_t^s \left(\frac{B_n^{lp}}{(T-u)^{lp/(l+1)}}+D_n^{l\bar{p}}+\abs{X_u}\right)du,
\end{eqnarray*}
with $p=1/(1-lp_g/(l+1))$ and $\bar{p}>1$. Since we assume that $lp_g <1$, then $lp/(l+1)<1$ and $\int_t^s \frac{B_n^{lp}}{(T-u)^{lp/(l+1)}}du \leqslant CB_n^{lp}$. Finally, we obtain
$$\int_t^s \abs{Z_u^M}^ldu \leqslant CA_n^l+CB_n^{lp}+CD_n^{l\bar{p}}+C\int_t^s\abs{X_u}du,$$
and 
$$\abs{X_s} \leqslant \abs{X_t}+C+C\int_t^s\abs{X_u}du+\sup_{t \leqslant r \leqslant T} \abs{\int_t^r \sigma(u) d\tilde{W}_u}+CA_n^l+CB_n^{lp}+CD_n^{l\bar{p}}.$$
Gronwall's lemma gives us
\begin{eqnarray*}
 \abs{X_s} &\leqslant& C\left(1+\sup_{t \leqslant r \leqslant T} \abs{\int_t^r \sigma(u) d\tilde{W}_u}+A_n^l+B_n^{lp}+D_n^{l\bar{p}}+\abs{X_t}\right)
\end{eqnarray*}
that implies
\begin{eqnarray}
\label{inegalite3 Z}
\mathbb{E}^{\mathbb{Q}^M}_t \left[\abs{X_s}^{c}\right] &\leqslant& C\left(1+A_n^{cl}+B_n^{clp}+D_n^{cl\bar{p}}+\abs{X_t}^c\right).
\end{eqnarray}
By putting (\ref{inegalite3 Z}) in (\ref{inegalite2 Z}) and (\ref{inegalite1 Z}), we obtain
\begin{eqnarray*}
 \abs{Z_t^M}^{l+1}(T-t) &\leqslant& C\left(1+\mathbb{E}^{\mathbb{Q}^M}_t \left[\abs{X_T}^{p_g}\right]+\int_t^T \mathbb{E}^{\mathbb{Q}^M}_t \left[\abs{X_s}^{p_g\vee (r_f+1)}\right]ds + (T-t)\abs{X_t}^{(l+1)r_f}\right)\\
 &\leqslant& C\left( 1+A_n^{(l+1)al}+B_n^{(l+1)alp}+D_n^{(l+1)al\bar{p}}+\abs{X_t}^{p_g}+(T-t)\abs{X_t}^{r_f+1} \right),
\end{eqnarray*}
with $a=(p_g\vee (r_f+1))/(l+1)$ and $C$ that does not depend on $M$ and constants that appear in assumption (TC.1). So, we easily see that we can take
$$A_{n+1}=C(1+A_n^{al}+B_n^{alp}+D_n^{al\bar{p}}), \quad B_{n+1}=C, \quad D_{n+1}=C,$$
and then the result is proved.

When $f$ and $g$ are not differentiable we can prove the result by a standard approximation and stability results for BSDEs with linear growth.
\cqfd

Since the estimate on $Z$ given by Proposition~\ref{prop estimation temporelle Z} does not depend on constants that appear in assumption (TC.1), we can use it to show an existence result for superquadratic BSDEs with a quite general terminal condition.

\begin{thm}
\label{theoreme final}
 Let assume that (F.1), (F.2), (B.1), (B.2)(b), (B.3) and (TC.2) hold. We also assume that $0 \leqslant p_gl <1$, then there exists a solution $(Y,Z)$ to the BSDE (\ref{EDSR}) such that $(Y,Z) \in \mathcal{S}^2 \times \mathcal{M}^2$. Moreover, we have for all $t \in [0,T[$,
\begin{equation}
\label{estimee deterministe 2 pour Z}
\abs{Z_t} \leqslant \frac{C(1+\abs{X_t}^{p_g/(l+1)})}{(T-t)^{1/(l+1)}}+C\abs{X_t}^{\frac{r_f+1}{l+1}},
\end{equation}
and, if we assume that (B.2)(c) holds,
\begin{equation*}
\label{moment Mp pour Z}
\E\left[\int_0^T\abs{Z_s}^{l+1}ds\right] <+\infty.
\end{equation*}
\end{thm}
\paragraph*{Proof of Theorem \ref{theoreme final}}
The proof is based on the proof of Proposition 4.3 in \cite{Delbaen-Hu-Bao-09}. For each integer $n\geqslant 0$, we construct the sup-convolution of $g$ defined by
$$g_n(x):=\sup_{u \in \R^d} \set{g(u)-n\abs{x-u}}.$$
Let us recall some well-known facts about sup-convolution:
\begin{lem}
 For $n\geqslant n_0$ with $n_0$ big enough, we have,
\begin{itemize}
 \item $g_n$ is well defined,
 \item (TC.1) holds for $g_n$ with $r_g=0$,
 \item (TC.2) holds for $g_n$ with same constants $C$ and $\bar{\alpha}$ than for $g$ (they do not depend on $n$),
 \item $(g_n)_n$ is decreasing,
 \item $(g_n)_n$ converges pointwise to $g$.
\end{itemize}
\end{lem}
Since (TC.1) holds, we can consider $(Y^n,Z^n)$ the solution given by Proposition \ref{existence unicite localement lipschitz}. It follows from Propositions \ref{comparison result} and \ref{estimation Y} that, for all $n \geqslant n_0$,
 \begin{equation}
 \label{croissance Yn}
 -C(1+\abs{X_t}^{p_g}+(T-t)\abs{X_t}^{r_f+1}) \leqslant Y^{n+1}_t \leqslant Y^n_t \leqslant Y^{n_0}_t \leqslant C(1+\abs{X_t}^{p_g}+(T-t)\abs{X_t}^{r_f+1}),
 \end{equation}
with $C$ that does not depend on $n$: indeed, the constant in Proposition \ref{estimation Y} just depends on the growth of the terminal condition and here the growth of $g_n$ can be chosen independently of $n$ (see previous lemma). So $(Y_n)_n$ converges almost surely and we can define 
$$Y=\lim_{n \rightarrow +\infty} Y^n.$$
Passing to the limit into (\ref{croissance Yn}), we obtain that the estimate of Proposition \ref{estimation Y} stays true for $Y$. Now the aim is to show that $(Z_n)_n$ converges in the good space. For any $T' \in]0,T[$, $(Y^n,Z^n)$ satisfies
\begin{equation}
 \label{etoile5}
Y_t^n=Y_{T'}^n+\int_t^{T'} f(s,X_s,Y_s^n,Z_s^n)ds-\int_t^{T'} Z_s^ndW_s, \quad 0 \leqslant t \leqslant T'.
\end{equation}
Let us denote $\delta Y^{n,m}:=Y^n-Y^m$ and $\delta Z^{n,m}:=Z^n-Z^m$. The classical linearization method gives us that $(\delta Y^{n,m}, \delta Z^{n,m})$ is the solution of BSDE
$$\delta Y^{n,m}_t = \delta Y^{n,m}_{T'}+\int_t^{T'} U_s^{n,m} \delta Y^{n,m}_s+V_s^{n,m} \delta Z^{n,m}_s ds -\int_t^{T'} \delta Z_s^{n,m} dW_s,$$
where $\abs{U^{n,m}} \leqslant C$ and, by using estimates of Proposition \ref{prop estimation temporelle Z},
\begin{equation}
\label{etoile1}
\abs{V^{n,m}} \leqslant C(1+\abs{Z^n}^l+\abs{Z^m}^l) \leqslant C(1+\abs{X}^{p}),
\end{equation}
with $p<1$ and $C$ that depends on $T'$ but does not depend on $n$ and $m$. Since $p<1$, Novikov's condition is fulfilled and we can apply Girsanov's theorem: there exists a probability $\Q^{n,m}$ such that $d\tilde{W}_t:=dW_t-V_t^{n,m}dt$ is a Brownian motion under this probability. By classical transformations, we have that $(\delta Y^{n,m}, \delta Z^{n,m})$ is the solution of the BSDE
$$\delta Y^{n,m}_t = \delta Y^{n,m}_{T'}e^{\int_t^{T'}U_s^{n,m}ds}-\int_t^{T'} e^{\int_t^{s}U_u^{n,m}du}\delta Z_s^{n,m} d\tilde{W}_s.$$
Since $U^{n,m}$ is bounded, classical estimates on BSDEs give us (see e.g. \cite{ElKaroui-Peng-Quenez-97})
\begin{equation}
 \label{etoile2}
\E^{\Q^{n,m}} \left[\left( \int_0^{T'} \abs{\delta Z^{n,m}_s}^2ds\right)^2 \right] \leqslant C\E^{\Q^{n,m}} \left[ \abs{\delta Y^{n,m}_{T'}}^4 \right].
\end{equation}
Now, we would like to have the same type of estimate than (\ref{etoile2}), but with the classical expectation instead of $\mathbb{E}^{\Q^{n,m}}$. To do so, we define the exponential martingale
$$\mathcal{E}^{n,m}_{T'}:= \exp \left( \int_0^{T'} V_s^{n,m} dW_s-\frac{1}{2}\int_0^{T'} \abs{V_s^{n,m}}^2ds \right).$$
Then, for all $p \in \R$, 
\begin{equation}
 \label{etoile3et4}
\E\left[ (\mathcal{E}^{n,m}_{T'})^p \right] <C_p,
\end{equation}
with $C_p$ that does not depend on $n$ and $m$: indeed, by applying (\ref{etoile1}) and Gronwall's lemma we have
\begin{eqnarray*}
 \E \left[e^{p \int_0^{T'} V_s^{n,m} dW_s-\frac{p}{2}\int_0^{T'} \abs{V_s^{n,m}}^2ds}\right] &=& \E \left[e^{\frac{1}{2}\left( \int_0^{T'} 2pV_s^{n,m} dW_s-\frac{1}{2}\int_0^{T'} \abs{2pV_s^{n,m}}^2ds\right) +\left(p^2-\frac{p}{2}\right) \int_0^{T'} \abs{V_s^{n,m}}^2ds}\right]\\
&\leqslant& \E \left[e^{\int_0^{T'} 2pV_s^{n,m} dW_s-\frac{1}{2}\int_0^{T'} \abs{2pV_s^{n,m}}^2ds}\right]^{1/2}\E\left[ e^{\left(2p^2-p\right) \int_0^{T'} \abs{V_s^{n,m}}^2ds}\right]^{1/2}\\
&\leqslant& \E\left[ e^{C\abs{2p^2-p}\left(1+\sup_{0 \leqslant s \leqslant T} \abs{X_s}^{2p}\right)}\right]^{1/2}\\
& < & + \infty,
\end{eqnarray*}
because $2p<2$. By applying Cauchy Schwarz inequality and by using (\ref{etoile3et4}) and (\ref{etoile2}), we obtain
\begin{eqnarray*}
 \E \left[\int_0^{T'} \abs{\delta Z^{n,m}_s}^2ds \right] &=& \E \left[(\mathcal{E}^{n,m}_{T'})^{-1/2}(\mathcal{E}^{n,m}_{T'})^{1/2}\int_0^{T'} \abs{\delta Z^{n,m}_s}^2ds\right]\\
&\leqslant& \E\left[ (\mathcal{E}^{n,m}_{T'})^{-1} \right]^{1/2}\E^{\Q^{n,m}} \left[\left( \int_0^{T'} \abs{\delta Z^{n,m}_s}^2ds\right)^2 \right]^{1/2}\\
&\leqslant& C\E^{\Q^{n,m}} \left[ \abs{\delta Y^{n,m}_{T'}}^4 \right]^{1/2}\\
&\leqslant& C\E\left[ (\mathcal{E}^{n,m}_{T'})^{2} \right]^{1/2}\E \left[ \abs{\delta Y^{n,m}_{T'}}^8 \right]^{1/4}\\
&\leqslant& C\E \left[ \abs{\delta Y^{n,m}_{T'}}^8 \right]^{1/4} \xrightarrow{n,m \to 0} 0.
\end{eqnarray*}
Since $\mathcal{M}^2$ is a Banach space, we can define
$$Z=\lim_{n \to +\infty} Z^n, \quad d\P \times dt\textrm{-a.e.} .$$
If we apply Proposition \ref{estimation Z}, we have that $\norm{Z^n}_{\mathcal{M}^2}<C$ with a constant $C$ that does not depend on $n$. So, Fatou's lemma gives us that $Z \in \mathcal{M}^2$. Moreover, the estimate on $Z^n$ given by Proposition \ref{prop estimation temporelle Z} stays true for $Z$ and, if we assume that (B.2)(c) holds, then Proposition \ref{estimation Z} gives us that
$$\E \left[ \int_0^T \abs{Z_s^n}^{l+1}ds\right]<C$$
with a constant $C$ that does not depend on $n$ and so 
$$\E \left[ \int_0^T \abs{Z_s}^{l+1}ds\right]<C.$$

Finally, by passing to the limit when $n \to +\infty$ in (\ref{etoile5}) and by using the dominated convergence theorem, we obtain that for any fixed $T' \in [0,T[$, $(Y,Z)$ satisfies 
\begin{equation}
 \label{etoile6}
Y_t=Y_{T'}+\int_t^{T'} f(s,X_s,Y_s,Z_s)ds-\int_t^{T'} Z_sdW_s, \quad 0 \leqslant t \leqslant T'.
\end{equation}
To conclude, we just have to prove that we can pass to the limit when $T' \to T$ in (\ref{etoile6}). Let us show that $Y_{T'} \xrightarrow{T' \to T} g(X_T)$ a.s.. Firstly, we have
$$\overline{\lim}_{s \to T} Y_s \leqslant \overline{\lim}_{s \to T} Y_s^n =g_n(X_T) \textrm{ a.s.} \quad \textrm{ for any } n \geqslant n_0, $$
which implies $\overline{\lim}_{s \to T} Y_s \leqslant g(X_T)$, a.s.. On the other hand, we use assumption (B.2)(b) and we apply Propositions \ref{estimation Y} and \ref{prop estimation temporelle Z} to deduce that, a.s.,
\begin{eqnarray*}
 Y_t^n &=& g_n(X_T)+\int_t^T f(s,X_s,Y_s^n,Z_s^n)ds -\int_t^T Z_s^n dW_s\\
&\geqslant& g_n(X_T)-C\int_t^T 1+\abs{X_s}^{r_f+1}+\abs{Y_s^n}+\abs{Z_s^n}^{\eta}ds-\int_t^T Z_s^n dW_s\\
&\geqslant& \E_t \left[ g_n(X_T) -C\int_t^T 1+\abs{X_s}^{(r_f+1) \vee p_g }+\frac{1+\abs{X_s}^{\eta p_g/(l+1)}}{(T-s)^{\eta/(l+1)}}ds\right]\\
&\geqslant& \E_t \left[ g_n(X_T)\right] -C(T-t)(1+\abs{X_t}^{(r_f+1) \vee p_g })-C(T-t)^{1-\eta/(l+1)}(1+\abs{X_t}^{\eta p_g/(l+1)}),
\end{eqnarray*}
and
$$Y_t=\lim_{n \to +\infty} Y_t^n \geqslant \E_t\left[ g(X_T)\right] -C(T-t)(1+\abs{X_t}^{(r_f+1) \vee p_g })-C(T-t)^{1-\eta/(l+1)}(1+\abs{X_t}^{\eta p_g/(l+1)}),$$
which implies
$$\underline{\lim}_{t \to T} Y_t \geqslant \underline{\lim}_{t \to T} \E_t\left[ g(X_T)\right] = g(X_T).$$
Hence, $\lim_{t \to T} Y_t=g(X_T)$ a.s. .

Now, let us come back to BSDE (\ref{etoile6}). Since we have
$$\int_t^T \abs{f(s,X_s,Y_s,Z_s)}ds \leqslant \int_t^T C(1+\abs{X_s}^{r_f+1}+\abs{Y_s}+\abs{Z_s}^{l+1})ds<+\infty \textrm{ a.s.},$$
then
$$\int_t^{T'} f(s,X_s,Y_s,Z_s)ds \xrightarrow{T' \to T} \int_t^{T} f(s,X_s,Y_s,Z_s)ds<+\infty\textrm{ a.s.}.$$
Finally, passing to the limit when $T' \to T$ in (\ref{etoile6}), we conclude that $(Y,Z)$ is a solution to BSDE (\ref{EDSR}). 
\cqfd

\begin{rem}
 The function $z \mapsto C\abs{z}^{l+1}+h(\abs{z}^{l+1-\eta})$ with $C>0$, $0<\eta\leqslant l+1$ and $h$ a differentiable function with a bounded derivative is an example of generator such that (B.1), (B.2)(b) and (B.3) hold.
\end{rem}

\begin{rem}
The estimate 
$$\abs{Z_t} \leqslant \frac{C(1+\abs{X_t}^{p_g})}{\sqrt{T-t}}+C\abs{X_t}^{r_f+1}$$
is already known in the Lipschitz framework as a consequence of the Bismut-Elworthy formula (see e.g. \cite{Fuhrman-Tessitore-02}). For the superquadratic case, the same estimate was obtained when $p_g=0$ and $f$ does not depend on $x$ and $y$ in \cite{Delbaen-Hu-Bao-09} (see also \cite{Richou-11} for the quadratic case). In \cite{Delbaen-Hu-Bao-09}, Remark 4.4. gives the same type of estimate than (\ref{estimee deterministe 2 pour Z}) for the example $f(z)=\abs{z}^{l}$. This result was already obtained by Gilding et al. in \cite{Gilding-Guedda-Kersner-03} using Bernstein's technique when $f(z)=\abs{z}^{l}$, $b=0$ and $\sigma$ is the identity.
\end{rem}

\begin{rem}
 In this article, estimate (\ref{estimee deterministe 2 pour Z}) for the process $Z$ allows us to obtain an existence result. But this type of deterministic bound is also interesting for numerical approximation of BSDEs (see e.g. \cite{Richou-11}) or for studying stochastic optimal control problems in infinite dimension (see e.g. \cite{Masiero-10}).
\end{rem}

\appendix
\section{Appendix}
\subsection{Proof of Lemma \ref{Qm sous martingale}}
Let us set
$$F_t^M := e^{\int_0^t\nabla_y f_M(s,X_s,Y_s^M,Z_s^M)ds}\nabla Y_t^M+\int_0^t e^{\int_0^s\nabla_y f_M(u,X_u,Y_u^M,Z_u^M)du}\nabla_x f_M(s,X_s,Y_s^M,Z_s^M)\nabla X_s ds,$$
and 
$$\tilde{F}_t^M := e^{\lambda t} F_t^M (\nabla X_t)^{-1}.$$
Since $d\nabla X_t=\nabla b(t,X_t)\nabla X_t dt$, then $d(\nabla X_t)^{-1} = - (\nabla X_t)^{-1}\nabla b(t,X_t)dt$ and thanks to Itô's formula,
$$d\tilde{Z}_t^M=dF_t^M(\nabla X_t)^{-1} \sigma(t)-F_t^M (\nabla X_t)^{-1}\nabla b(t,X_t)\sigma(t)dt+F_t^M(\nabla X_t)^{-1}\sigma'(t)dt,$$
and
$$d(e^{\lambda t}\tilde{Z}_t^M)=\tilde{F}_t^M(\lambda Id-\nabla b(t,X_t))\sigma(t)dt+\tilde{F}_t^M\sigma'(t)dt+e^{\lambda t} dF_t^M(\nabla X_t)^{-1} \sigma(t).$$
Finally,
$$d\abs{e^{\lambda t} \tilde{Z}_t^M}^2=d\langle N \rangle_t+2\left[\lambda\abs{\tilde{F}_t^M\sigma(t)}^2-\tilde{F}_t^M\sigma(t)[\tr{\sigma(t)}\tr{\nabla b(t,X_t)}-\tr{\sigma'(t)}]\tr{\tilde{F}_t^M}\right]dt+dN_t^*,$$
with $N_t:=\int_0^t e^{\lambda s} dF_s^M(\nabla X_s)^{-1}\sigma(s)$ and $N_t^*$ a $\mathbb{Q}^M$-martingale. Thanks to the assumption (F.2) we are able to conclude that $\abs{e^{\lambda t} \tilde{Z}_t^M}^2$ is a $\mathbb{Q}^M$-submartingale.
\cqfd

\def\cprime{$'$}

\end{document}